\documentclass[10pt,a4paper]{article}
\usepackage{amssymb,amsmath}
\usepackage{t1enc}
\newtheorem{tw}{Theorem}[section]
\newtheorem{df}{Definition}[section] 
\newtheorem{st}{Proposition}[section]
\newtheorem{prz}{Example}[section]
\newtheorem{wn}{Corollary}[section]
\def \bw {\begin{wn}}
\def \ew {\end{wn}}
\def \bprz {\begin{prz}}
\def \eprz {\end{prz}}
\def \bst {\begin{st}}
\def \est {\end{st}}
\def \bt {\begin{tw}} 
\def \et {\end{tw}} 
\def \bd {\begin{df}}
\def \ed {\end{df}}
\def \brt {\begin{flushright}}
\def \ert {\end{flushright}}
\def \bre {\blacksquare} 
\def \it {\indent} 
 
\def \sbt {\subset}

\def \ra {\rightarrow} 
\def \Ra {\Rightarrow} 
\def \sm {\setminus} 
\def \es {\emptyset}
\def \lo {\leadsto} 
\def \mc {\mathcal} 
\def \tf {\textbf} 
\def \oe {\overline} 
\def \li {\limits}
\def \iy {\infty}
\def \ve {\varepsilon}
\def \la {\lambda}
\def \aa {\alpha}

\def \da {\delta}

\def \z {\zeta}
\def \g {\gamma}
\def \xrw {\xrightarrow}
\def \bm {\bibitem}
\def \ie {\itshape}
\def \c {\cite}
\def \mb {\mathbb} 

\def\blist{\begin{list}{}{\setlength{\partopsep}{0pt}\setlength{\itemsep}{0pt}%
\setlength{\topsep}{0pt}\setlength{\parskip}{0pt}}}
\def\elist{\end{list}}
\def \cn {\colon} 

\title{Continuous version of the Choquet Integral Representation
Theorem}
\author{Piotr Pucha{\l}a}
\begin{document}
\setcounter{page}{1}
\maketitle

\begin{abstract}
Let $E$ be a locally convex topological Hausdorff space, 
$K$ --- its nonempty compact convex subset, $\mu$ --- a regular,
probability Borel measure on $E$ and $\g >0$.
We say that the measure $\mu$ $\g$ ---
represents point $x\in K$,  
if for any $f\in E^{\ast}$  an inequality
$\sup\li_{\| f\|\leq1}| f(x)-\int\li_{K}fd\mu | <\g$ holds.
In this paper the continuous version of the Choquet theorem is
proved. Namely, it is shown that if $P$ is continuous multivalued
mapping from a metric space $T$ into the space of nonempty, bounded
convex subsets of a Banach space $X$, then there exists a 
$\text{\text{weak}}^{\ast}$ continuous family $(\mu_{t})$ of regular Borel
probability measures on $X$ $\g$ -- representing points in
$P(t)$. The
two cases are considered: in the first one the values of $P$ are
compact while in the second --- closed. For this purpose it is shown
(using geometrical tools)
that the mapping $t\ra\text{\text{ext}}P(t)$ is lower semicontinuous.
The continuous versions of the Krein -- Milman theorem
are obtained as corollaries.
\end{abstract}



\section{Introduction}
\it\it The classical Minkowski -- Carath\'eodory representation
theorem states that each point of a compact convex set $K$ in
$\mb{R}^{n}$ can be written as a convex combination of at most $n+1$
extreme points of $K$. This theorem was ge\-ne\-ra\-li\-zed by
G. Choquet 
(\c{choq1}) who proved that  each $k$ point of a compact, convex and 
metrizable subset $K$ of a locally convex Hausdorff topological space
$X$ is a barycenter of a regular Borel probability measure $\mu_{k}$
on $X$, 
supported by extreme points points of $K$, i.e. that the equality 
 \[
    f(k)=\int\li_{K}fd\mu_{k}
 \]
holds for any $f\in X^{\ast}$ with 
$\mu_{k}(\text{\text{ext}}K)=1$, where $\text{\text{ext}}K$ stands for
the set of extreme points of $K$.

E. Bishop and K. de Leeuw (\c{bis}) removed the metrizability
assumption.
G. A. Edgar (\c{ed}) proved the noncompact version of
Choquet theorem. His result stated that the thesis of 
Choquet -- Bishop -- de Leeuw remained true for $K$ being nonempty
bounded, closed, convex and separable subset of a~Banach space $X$
having Radon -- Nikod\'ym property (RNP for short). In the paper
\c{ed1} he improved his result by removing separability
condition. P.Mankiewicz (\c{ma}) modified it by introducing
``separable extremal ordering'', more natural and easy to use partial
order than the one introduced by Edgar (see \c{bour}, p.174).

The purpose of this paper is to show that an analogue of the Choquet
theorem holds for ''moving'' sets. These sets are values of a
multivalued mapping, defined on  a metric space $T$ into
suitable subsets of a certain Banach space $X$. We shall  consider
two cases.

In the first one 
 multivalued mapping $P\cn T\lo X$ is continuous with
compact convex values in a separable Banach space $X$. In the second
 case $X$ is 
separable, reflexive Banach space while $P\cn T\lo X$ is continuous with
bounded closed convex values (recall that ref\-le\-xive Banach space has
RNP). By the celebrated Michael theorem there exists a continuous
function $p\cn T\ra X$, called a continuous selection of $P$, with the property  
that $p(t)\in P(t)$ for all $t\in T$. It will be shown that each such point
$p(t)$ is, for given selection $p(\cdot)$, an ,,almost barycenter'' of
the regular Borel probability measure $\mu_{t}$ on $X$ such that
$\mu (\text{\text{ext}}P(t))=1$. In other words, for any continuous
selection  $p(t)$ of the multifunction $P(t)$ there exists a
continuous (in the weak* topology) family of measures 
$(\mu_{t})_{t\in T}$ ,,almost representing'' points $p(t)$. In both
cases the fact that multifunction
 \[
t\ra\normalfont{\text{ext}}P(t)
 \]
is lower semicontinuous is crucial for the main result. As obvious
corollaries we obtain continuous versions of the Krein -- Milman
theorem. 

All the necessary information about multifunctions can be found in 
\c{hu1}, for Choquet theorem see \c{alf}, \c{ph} and \c{bour}
(noncompact case), Banach spaces with Radon -- Nikodym property are
subject of the classics \c{bour}, \c{du}, properties of measures on
metric spaces are investigated in \c{ali} (where one can also find a
chapter devoted to multivalued mappings).

\section{Preliminaries}
\it\it In this section we state several definitions and facts needed.
\bd
Let $X$ and $Y$ be topological spaces and $P\cn X\lo Y\sm\es$ --- 
a~set -- valued map. We say that
\begin{description}
\item[(a)] $P$ is lower semicontinuous (lsc) iff the set
 \[
P^{-}(U):=\{ x\in X:P(x)\cap U\neq\es\}
 \]
is open whenever $U\sbt Y$ is open;
\item[(b)] $P$ is upper semicontinuous (usc) iff the set 
 \[
P^{-}(V):=\{ x\in X:P(x)\sbt U\}
 \]
is closed whenever $V\sbt Y$ is closed;
\item[(c)] $P$ is continuous iff it is both lower -- and upper
semicontinuous. 
\end{description}
\ed
\bt (Michael) Let $X$ be a paracompact space, $Y$ -- a Banach space
and $P\cn X\lo Y$ -- lower semicontinuous multifunction with convex
values. Then:
\begin{description}
\item[(a)] for any $\ve >0$ there exists a continuous function
$p_{\ve}\cn X\ra Y$such that $d(p_{\ve}(x),P(x))<\ve$ 
$\forall x\in X$; this function is called an $\ve$ -- selection of
$P$;
\item[(b)] if in addition the values of $P$ are closed, then there
exists a continuous function $p\cn X\ra Y$ such that $p(x)\in P(x)$;
this function is called a~continuous selection of $P$.
\end{description}
\et

Now denote by $X$ a locally convex topological Hausdorff space and let
$K$ be its compact convex subset. If $\mu$ is a regular, Borel
probability measure on $X$, we say that it is supported by the set
$A\sbt X$ (not necessarily closed) if $\mu (A)=1$.
\bd 
For such $X$, $K$ and $\mu$ be as above. We say that:
\begin{description}
\item[(a)] measure $\mu$ represents point $x\in K$, if for all 
$f\in X^{\ast}$ we have
 \[
f(x)=\int\li_{K}fd\mu.
 \]
This point, denoted by $r(\mu)$, is called the barycenter of $\mu$;
\item[(b)] measure $\mu$ $\g$ --- represents point $x\in K$, $\g >0$,
if for any $f\in X^{\ast}$  an inequality
 \[ 
\sup\li_{\| f\|\leq1}\left| f(x)-\int\li_{K}fd\mu\right| <\g .
 \]
holds.
\end{description}
\ed
\bt (Choquet) Let $X$, $K$ and $\mu$ be as above and assume
additionally that $K$ is metrizable. Then for any $x\in K$ there
exists a regular Borel probability measure $\mu_{x}$ representing point
$x$ and such that $\mu_{x}(\normalfont{\text{ext}}K)=1$.
\et
Recall that if the set $\normalfont{\text{ext}}K$ is closed then the
Choquet theorem is equivalent to the Krein -- Milman theorem.
\bt (Edgar -- Mankiewicz noncompact version of the Choquet theorem)
Let $K$ be a (nonempty) closed bounded convex subset of a Banach space
$X$ and suppose that $K$ has RNP. Then every point of $K$ is a
barycenter of a~regular Borel probability measure on $K$ supported by
the set $\normalfont{\text{ext}}K$.
\et
\section{Compact case}
\it\it In this section we deal with multifunction $P$ from $T$ into
compact convex sets of $X$.

We first establish the lower semicontinuity of map with values in the
set of extreme points of the compact convex set. Recall that exposed
point of a compact convex subset of a Banach space is a strongly
exposed point of this subset.
\bst
Let $T$ be a metric space, $X$ --- a Banach space, 
$P\cn T\lo X$ --- a continuous multifunction with compact convex
values. Then the multifunction
 \[
t\ra\normalfont{\text{ext}}P(t)
 \]
is lower semicontinuous.
\est
\tf{Proof.} Let $(t_{n})$ be a sequence in $T$, convergent to a
point $t_{0}\in T$. It is enough to show that for each such sequence
and any $a_{0}\in\text{\text{ext}}P(t_{0})$ there exists sequence
$(a_{n})$, $n\in\mb{N}$, such that $a_{n}\in\text{\text{ext}}P(t_{n})$
and $a_{n}\xrw[n\ra\iy]{}a_{0}$.

Let $e_{0}$ be any exposed (in fact strongly exposed) point of
$P(t_{0})$. Then there exists functional $f_{0}\in X^{\ast}$, with unit
norm, strongly exposing $e_{0}$. The lower semicontinuity of $P$
yields existence of a sequence $(x_{n})\sbt X$, convergent to
$e_{0}$, with $x_{n}\in P(t_{n})$. Fix a number $\g >0$ and define
the slice
 \[
R_{\g}(t_{n}):=\{x\in P(t_{n}):f_{0}(x)>c(f_{0},P(t_{n}))-\g\},
 \]
where $c(\cdot,A)$ stands for the support function of the set $A$. It
turns out that there exists $n_{0}\in\mb{N}$ such that for every
$n\geq n_{0}$ the intersection 
$R_{\g}(t_{n})\cap\text{\text{ext}}P(t_{n})$ is nonempty. Suppose not. 
Then for each $n_{0}\in\mb{N}$ there exists $n\geq n_{0}$ for which
$R_{\g}(t_{n})\cap\text{\text{ext}}P(t_{n})=\es$, i.e.
$\text{\text{ext}}P(t_{n})\sbt X\sm R_{\g}(t_{n})$. This implies the
existence of a subsequence $n_{k}\xrw[k\ra\iy]{}\iy$ having
property that for each $e\in\text{\text{ext}}P(t_{n_{k}})$ 
 \[
f_{0}(e)<c(f_{0},P(t_{n_{k}}))-\g.
 \]
By the Krein -- Milman theorem the set $P(t_{n_{k}})$ coincides with the
closed convex hull of its extreme points, so in particular we have
 \[
f_{0}(x_{n_{k}})\leq c(f_{0},P(t_{n_{k}}))-\g,\;\;\; k=1,2,...
 \]
Passing to the limit we obtain the inequality
 \[
f_{0}(e_{0})\leq c(f_{0},P(t_{0}))-\g,
 \]
contradicting the fact, that $e_{0}$ is strongly exposed.

Now let $\g =\tfrac{1}{m}$, $m=1,2,...$, and consider slices
$R_{\tfrac{1}{m}}(\cdot)$. Then for each $m$ there exists $n_{m}$ such
that for $n\geq n_{m}$ we have
 \[
R_{\tfrac{1}{m}}(t_{n})\cap\normalfont{\text{ext}}P(t_{n})\neq\es.
 \]
We can assume that $n_{m}\leq n<n_{m+1}$. For such $n$ choose
$e_{n}\in R_{\tfrac{1}{m}}(t_{n})\cap\normalfont{\text{ext}}P(t_{n})$
obtaining the sequence $(e_{n})$ with the property
 \[
f_{0}(e_{n})\geq c(f_{0},P(t_{n}))-\tfrac{1}{m}
 \]
for $n_{m}\leq n<n_{m+1}$.

By the upper semicontinuity of $P$ there exists subsequence of
$(e_{n})$ (denoted also $(e_{n})$) convergent to a point 
$\oe{e}_{0}\in P(t_{0})$. As the values of $P$ are convex sets, we can
use the relationship between the Hausdorff distance $h(A,B)$ of the
sets $A$, $B$ and their support functions, so 
 \[
\sup\li_{\| f\|\leq 1}\{|c(f,P(t_{n}))-c(f,P(t_{0}))|\}=
h(P(t_{n}),P(t_{0})).
 \]
The compactness of $P(t)$, $t\in T$, yields the continuity of $P$ in
Hausdorff metric, which gives us 
$h(P(t_{n}),P(t_{0}))\xrw[n\ra\iy]{}0$, which in turn implies that 
for any $f\in X^{\ast}$ 
$c(f,P(t_{n}))\xrw[n\ra\iy]{}c(f,P(t_{0}))$. But for 
$n_{m}\leq n\leq n_{m+1}$ we have 
$f_{0}(e_{n})\geq c(f_{0},P(t_{n}))-\tfrac{1}{m}$. Taking into account
that $f_{0}(e_{n})$ converges to $f_{0}(\oe{e})$ we finally get
$e=\oe{e}_{0}$.

Now by Lindenstrauss -- Troyanski result (\c{lin}, \c{tr}, see also
\c{bour}) the set $P(t_{0})$ equals the closed convex hull of its
(strongly) exposed points and by Milman partial converse of the Krein
-- Milman theorem those points are dense in the set
$\text{\text{ext}}P(t_{0})$. 

We are now ready to construct a desired sequence of extreme points.
So let $a_{0}$ be any extreme point of $P(t_{0})$. Choose and fix
$n_{1}\in\mb{N}$ and $e_{0}^{n_{1}}\in\normalfont{\text{(st)exp}}P(t_{0})$.
There exists sequence $(b_{n}^{1})$ of extreme points of $P(t_{0})$
convergent to $e_{0}^{n_{1}}$. Then there exists $n_{2}>n_{1}$ such
that for $n\geq n_{2}$ we have
$\|e_{0}^{n_{1}}-b_{n}^{1}\|<\tfrac{1}{n_{1}}$. Now take 
$e_{0}^{n_{2}}\in\normalfont{\text{(st)exp}}P(t_{0})$ with
$\|e_{0}^{n_{2}}-a_{0}\|<\tfrac{1}{n_{2}}$ and sequence $(b_{n}^{2})$
of extreme points of $P(t_{n})$ convergent to $e_{0}^{n_{2}}$. Then
there exists such $n_{3}>n_{2}$ that for all $n\geq n_{3}$ we have the
inequality $\|e_{0}^{n_{2}}-b_{n}^{2}\|<\tfrac{1}{n_{1}}$. Repeating
this procedure we obtain sequences $(b_{n}^{i})$, 
$b_{n}^{i}\in\text{\text{ext}}P(t_{n})$. Setting $a_{n}:=b_{n}^{i}$ we
obtain the desired sequence.\mbox{}\hfill$\bre$

\tf{\textsc{Remark.}} Tolstonogov and Figonienko proved (under the
same assumptions) in \c{to} the lower semicontinuity of the map
$t\ra\text{\text{clext}}P(t)$, where ,,cl'' stands for
,,closure''. This result is equivalent to the above one, but the
method of the proof presented here is of geometrical nature, in
contrast to the topological methods they used. Incidentally, it seems
that considering the map $t\ra\text{\text{ext}}P(t)$, instead of
$t\ra\text{\text{clext}}P(t)$ is ,,in the spirit'' of the Choquet
theorem. 

Now let $T$ be a metric space, $X$ --- separable Banach space. By
$\mc{M}(X)$ we denote regular probability Borel measures on $X$.
We consider a continuous multifunction
 \[
P\cn T\lo X,
 \]
 with (nonempty) compact convex values. The Michael
selection theorem assures us of the existence of a continuous
selection $p$ of $P$. Define the set -- valued map
$L\cn T\lo\mc{M}(X)$ setting
 \[
L(t):=\biggl\{\mu\in\mc{M}(X):\mu
  (\normalfont{\text{ext}}P(t))=1,\;\;\; 
\sup\li_{\|f\|\leq 1}
\biggl|f(p(t))-\int\li_{P(t)}fd\mu\biggr|<\g ,\;\;\; 
f\in X^{\ast}\biggr\},
 \]
where $\g$ is fixed positive number. Choquet theorem guarantees that
$L(t)$ is nonempty for all $t\in T$. Choose and fix the continuous
selection $p$ of $P$.
\bst
Multifunction $L$ is lower semicontinuous.
\est
\tf{Proof.} It is enough to show that for any sequence $(t_{n})$ of
$T$, convergent to $t_{0}\in T$ and for any nonempty,
$\text{\text{weakly}}^{*}$ closed subset $F$ of $\mc{M}(X)$,
 an implication $L(t_{n})\sbt F\Ra L(t_{0})\sbt F$ holds.

Take any element $\mu_{0}$ of $L(t_{0})$.
The set of discrete measures on $\text{\text{clext}}P(t_{0})$
is dense in the set of measures supported by that set, so there exists
a sequence $(m_{k})$ of discrete measures (i.e. convex combinations of
Dirac measures), convergent
$\text{\text{weakly}}^{\ast}$ to $\mu_{0}$. This yields the existence of
such $k_{0}$ that for all $k\geq k_{0}$ measure $m_{k}$ 
$\g-\text{\text{represents}}$ point $p(t_{0})$. Each measure $m_{k}$, 
$k\geq k_{0}$, is of the form
 \[
m_{k}=\sum\li_{i=1}^{m}\la_{i}\da_{a_{i}},
 \]
where $\da_{(\cdot)}$ is Dirac measure,
$a_{i}\in\text{\text{ext}}P(t_{0})$ and $\la_{i}$ are coefficients of
the convex combination. As $\text{\text{ext}}P(\cdot)$ is lower
semicontinuous, for any $a_{i}$, $i=1,2,...,m$, there exists a sequence
$(b_{n}^{i})$, $b_{n}^{i}\in\text{\text{ext}}P(t_{n})$, convergent to
$a_{i}$. This means that for fixed $k$ the sequence $(\mu_{n}^{k})$ of
measures, $\mu_{n}^{k}=\sum\li_{i=1}^{m}\la_{i}\da_{b_{n}^{i}}$
converges to $m_{k}$. For $f\in X^{\ast}$ we also have
\begin{multline}
\biggl| f(p(t_{n})) -\int\li_{P(t_{n})}fd\mu_{n}^{k}\biggr|\leq\\\nonumber
\leq |f(p(t_{n})) - f(p(t_{0})) |+
\biggl| f(p(t_{0})) -\int\li_{P(t_{0})}fd\mu_{0}\biggr|+\\\nonumber
+\biggl|\int\li_{P(t_{0})}fd\mu_{0}
-\int\li_{P(t_{0})}fdm_{k}\biggr|+
\biggl|\int\li_{P(t_{0})}fdm_{k}
-\int\li_{P(t_{n})}fd\mu_{n}^{k}\biggr|.
\end{multline}
The first and the last terms on the right converge to zero. Since there
exists  $M>0$ such that
$\sup\{\|y\| :y\in P(t_{0})\}\leq M$, so for $x\in P(t_{0})$ we have
$|f(x)|\leq M|f(\tfrac{1}{M}x)|\leq M$. Taking into account that the
sets $\text{\text{supp}}m_{k}$ and $\text{\text{supp}}\mu_{0}$ are the
subsets of $P(t_{0})$ we obtain
\begin{equation}
 \sup\li_{\| f\|\leq 1}\biggl|\int\li_{P(t_{0})}fdm_{k}
-\int\li_{P(t_{0})}fd\mu_{0}\biggr|=
\sup\li_{\| f\|\leq 1}
\biggl|\int\li_{P(t_{0})}fd(m_{k}-\mu_{0})\biggr|\leq
M|(m_{k}-\mu_{0})(1)|\xrightarrow[k\ra\iy]\quad 0.
\end{equation}
For fixed $k\geq k_{0}$ and for all $n\geq n_{0}$ we then have
\[\biggl| f(p(t_{n})) -\int\li_{P(t_{n})}fd\mu_{n}^{k}\biggr|<\g. \]
By construction, measure $\mu_{n}^{k}$ is supported by the set of
extreme points of $P(t_{n})$, so for $n\geq n_{0}$ it
belongs to $L(t_{n})$ and thus to $F$. Passing to the limit gives
inclusion $\mu_{0}\in F$ resulting the lower semicontinuity
of the multifunction $L$.
\mbox{}\hfill$\bre$
\bw 
Multifunction
\[\normalfont{\text{cl}}L(t)=
 \biggl\{ \mu\in\mc{M}(X):\,
\mu (\normalfont\text{{clext}}P(t))=1,\,
\sup\li_{\|f\|\leq 1}\biggl| f(p(t)) -\int\li_{P(t)}fd\mu\biggr| <\g,
\;\; f\in X^{\ast }\biggr\} \]
is lower semicontinuous.
\ew
\bw
There exists a continuous selection
$\oe{l}\cn T\ra \mc{M}(X)$ of the multifunction 
$\normalfont{\text{cl}}L$.
\ew

The set of extreme points of even compact convex set need not be
closed, so in general the values of the multifunction $L$ are not
closed sets. In particular we cannot expect $L$ to have continuous
selections. However, there exists continuous approximate selections,
as stated in the next result.

By $\mc{H}$ we will denote the Hilbert cube, $\mc{M}(\mc{H})$ stands
for the regular probability Borel measures on $\mc{H}$; this is compact
separable metric space.
\bt (continuous version of Choquet Theorem)\\
Let $T$ be a metric space and $X$ --- separable Banach space. Consider
continuous multifunction $P\cn T\lo X$ with compact convex values and
denote by $p$ its continuous selection. Consider also multifunction
$L\cn T\lo\mc{M}(X)$:
\[ L(t):=\biggl\{\mu\in\mc{M}(X):
\mu (\normalfont{ext}P(t))=1,\,
 \sup\li_{\|f\|\leq 1}\biggl| f(p(t))-\int\li_{P(t)}fd\mu\biggr| <\g,\;\;
f\in X^{\ast}\biggr\}. \]
Then for any $\da >0$ $L$ admits a continuous
$\da -\text{\text{selection}}$ of
 $L$.
\et
\tf{Proof.} In what follows we adopt classical Michael method.
There exists (see \c{ali}, pp. 483 -- 485) continuous  mapping 
$\hat{\phi}\cn\mc{M}(X)\ra\mc{M}(\mc{H})$. The Polish space
 $\mc{M}(X)$ is metrizable by
 \[ d(\mu_{1},\mu_{2})=
 \rho_{\mc{H}}(\hat{\phi}(\mu_{1}),\hat{\phi}(\mu_{2}))= 
 \sum\li_{j=1}^{\iy}\frac{1}{2^{j}}
 |(\hat{\phi}(\mu_{1})-\hat{\phi}(\mu_{2}))(\z_{j})|, \]
where $\z_{j}$, $j=1,2,...$, are from the dense set in $C(\mc{H})$ --
the space of the continuous functions on the Hilbert cube. We start
with fixing: a dense set $(\z_{n})_{n=1}^{\iy}$ in $C(\mc{H})$,
 numbers $\da >0$ and $N\in\mb{N}$ with 
$\sum\li_{j=N+1}^{\iy}\frac{1}{2^{j}}<\frac{\da}{4}$, functions
$\z_{1},...,\z_{n}\in (\z_{n})_{n=1}^{\iy}$, point $t_{0}\in T$ and
measure $\mu_{o}\in L(t_{0})$.
The mapping $\mu\ra |(\hat{\phi}(\mu)-\hat{\phi}(\mu_{0}))(\z_{j})|$
is continuous for (any) $\mu_{0}$ and any fixed $j\in\mb{N}$, so the
 set 
\[ V(\mu_{0},\z_{1},\dots,\z_{N},\da):=
 \biggl\{\mu:\sum\li_{j=1}^{N}\frac{1}{2^{j}}
 \left|(\hat{\phi}(\mu)-\hat{\phi}(\mu_{0}))(\z_{j})\right|<\frac{\da}{2}\biggr\}
 \] 
is open. The lower semicontinuity of $L$ implies that 
\[ U(t_{0},\mu_{0}):=L^{-1}(V)=\{t\in T:L(t)\cap V\neq\es\} \]
is an open neighborhood of $t_{0}$ is nonempty and open. Using the
lower semicontinuity of $L$ again 
we obtain an open cover $\{U(t_{\aa},\mu_{\aa})\}_{\aa\in I}$ of $T$.
By $e_{\aa}(\cdot)$
denote locally finite partition of unity subordinated to this
covering. Our candidate for continuous $\da-\text{\text{selection}}$
$l_{\da}$ of $L$ is of the form
 \[ l_{\da}(t):=\sum\li_{\aa\in I}e_{\aa}(t)\mu_{\aa}. \]
Fix $t\in T$ and set $\{\aa :e_{\aa}(t)>0\}:=\{\aa_{1},...,\aa_{k}\}$;
then $t\in\text{\text{supp}}e_{\aa_{i}}\sbt
U(t_{\aa_{i}},\mu_{\aa_{i}})$, so the intersection of $L(t)$ with ball
of radius $\da$, centered in $\mu_{\aa_{i}}$, is nonempty. Now take
measure $\oe{\mu}_{i}$ -- the element of this intersection, and
observe that
$d(\oe{\mu}_{i},\mu_{\aa_{i}})<\da$. The point
$\oe{l}_{\da}(t):=\sum\li_{i=1}^{k}e_{\aa_{i}}(t)\oe{\mu}_{i}$ lies in
the convex set $L(t)\cap V$. The thesis of the theorem follows from
the following sequence of inequalities:
 \begin{align}
 d(l_{\da}(t),L(t))& \leq d(l_{\da}(t),\oe{l}_{\da}(t)) =
 \rho_{\mc{H}}(\hat{\phi}(l_{\da}(t)),\hat{\phi}(\oe{l}_{\da}(t)))\leq
 \nonumber\\[.5cm] 
 & \leq\sum\li_{j=1}^{\iy}\frac{1}{2^{j}}
 \left|\left(\sum\li_{\aa\in I}e_{\aa}(t)\hat{\phi}(\mu_{\aa})-
 \sum\li_{i=1}^{k}e_{\aa_{i}}(t)
 \hat{\phi}(\oe{\mu}_{i})\right)\left(\z_{j}\right)\right|\leq \nonumber\\
 & \leq\sum\li_{j=1}^{\iy}\frac{1}{2^{j}}\sum\li_{i=1}^{k}e_{\aa_{i}}(t)
 \left|(\hat{\phi}(\mu_{\aa_{i}})-\hat{\phi}(\oe{\mu}_{i}))(\z_{j})\right|\leq
 \nonumber \\
 &\leq \sum\li_{j=1}^{N}\frac{1}{2^{j}}
 \sum\li_{i=1}^{k}\left|(\hat{\phi}(\mu_{\aa_{i}})-
 \hat{\phi}(\oe{\mu}_{i}))(\z_{j})\right|
  +\sum\li_{j=N+1}^{\iy}\frac{1}{2^{j}}< \nonumber\\
 & < \sum\li_{i=1}^{k}e_{\aa_{i}}(t)\left(\sum\li_{j=1}^{N}
 \left|(\hat{\phi}(\mu_{\aa_{i}})-\hat{\phi}(\oe{\mu}_{i}))(\z_{j})\right|\right)+
 \frac{\da}{4}<\nonumber\\
 &<\frac{\da}{2}+\frac{\da}{4}<\da.\nonumber
 \end{align}
\mbox{}\hfill$\bre$
\bw
(continuous version of the Krein -- Milman theorem)\\
Let $T$, $X$, $P$ and $p$ be as above. 
Then for any $\g>0$ there exists continuous family of
measures $(\mu_{t})_{t\in T}$ on $X$, supported by the closure of the
set of extreme points of $P(t)$ and $\g -\text{\text{representing}}$
point $p(t)$.
\ew
\section{Noncompact case}
\it\it In this section we consider the noncompact case. We have to
impose more assumptions both on the Banach space $X$ and multifunction
$P$.

Recall that $x\in K$ is a denting point iff for each $\ve >0$ 
$x\not\in\text{\text{clcv}}(K\sm U_{\ve}(x))$, where $U_{\ve}(x)$
denotes the $\ve$ --- neighbourhood of $x$. 
\bt
Let $T$ be a metric space, $X$ --  reflexive Banach space
and consider multifunction $P\cn T\lo X$ fulfilling the following
conditions:
\begin{description}
\item[(a)] $P$ is continuous;
\item[(b)] for each $t\in T$ the set $P(t)$ is bounded, closed and
convex;
\item[(c)] each extreme point of the set $P(t)$ is its denting point.
\end{description}
Then the multifunction
 \[
t\ra\normalfont{\text{ext}}P(t)
 \]
is lower semicontinuous.
\et
\tf{Proof.} Exactly as in the compact case we construct the slice
$R_{\g}(\cdot)$ and show that it is nonempty, replacing '' Krein --
Milman theorem'' with ''Krein -- Milman property'', obtaining in the
same way the sequence $(e_{n})$ of extreme points belonging for 
$n_{n}\leq n\leq n_{m+1}$ both to the set $\text{\text{ext}}P(t_{n})$
and the slice $R_{\tfrac{1}{m}}(t_{n})$. By the upper semicontinuity
of $P$ this sequence is bounded and has a  subsequence
(denoted also by $(e_{n})$) convergent to the point 
$\oe{e}_{0}\in P(t_{0})$. The upper semicontinuity of $P$ implies its
Hausdorff upper semicontinuity, so we can write
\[
   \sup\{ c(f,P(t_{n}))-c(f,P(t_{0})):\|f\|\leq 1\}=
h^{*}(P(t_{n}),P(t_{0}))\xrw[n\ra\iy]{} 0,
 \]
where $h^{\ast}(A,B)=\sup\{d(a,B):a\in A\}$.
Thus $c(f,P(t_{n}))\xrw[n\ra\iy]{}c(f,P(t_{0}))$ and (as in the
compact case) we can conclude that $\oe{e}_{0}=e_{0}$.

We have thus costructed the sequence $(e_{n})$, 
$e_{n}\in\text{\text{ext}}P(t_{n})$, weakly convergent to the point 
$e_{0}\in\text{\text{stexp}}P(t_{0})$. Moreover, we have\\
$d(e_{n},P(t_{0}))\leq h^{\ast}(P(t_{n}),P(t_{0}))
\xrw[n\ra\iy]{}0$, so there exists sequence $(b_{n})$ of elements of
$P(t_{0})$ with $\|b_{n}-e_{n}\|\xrw[n\ra\iy]{}0$. This yields weak
convergence of $(b_{n})$ to $e_{0}$ which in turn implies (as $e_{0}$
is strongly exposed), that $(b_{n})$ converges to $e_{0}$ in norm.
Thus $\|e_{n}-e_{0}\|\xrw[n\ra\iy]{}0$.

It turns out that the set of strongly exposed points of $P(t)$, $t\in T$,  
is dense in the set $\text{\text{ext}}P(t)$. Indeed, suppose that this
is not the case and consider slice  of
$P(t)$ with norm diameter $\ve$ containing some
$e\in\text{\text{ext}}P(t)$. None of the strongly 
exposed points of $P(t)$ belongs to the slice, so we have 
$P(t)=\text{\text{clcvstexp}}P(t)$. This contradicts the fact that $e$
is the denting point.

The rest of the proof proceeds as in the compact case.
\mbox{}\hfill$\bre$

Now we are able to reformulate continuous version of the Choquet
theorem and its corollary.

\bt ( continuous version of the noncompact Choquet theorem)\\
Let $T$ be a metric space, $X$ -- separable, reflexive Banach space
and consider multifunction $P\cn T\lo X$ fulfilling the following
conditions:
\begin{description}
\item[(a)] $P$ is continuous;
\item[(b)] for each $t\in T$ the set $P(t)$ is bounded, closed and
convex;
\item[(c)] each extreme point of the set $P(t)$ is its denting point.
\end{description}
Denote by $p$ the continuous selection of $P$ and define the set --
valued map 
$L\cn T\lo\mc{M}(X)$:
\[ L(t):=\biggl\{\mu\in\mc{M}(X):
\mu (\normalfont{ext}P(t))=1,\,
 \sup\li_{\|f\|\leq 1}\biggl| f(p(t))-\int\li_{P(t)}fd\mu\biggr| <\g,\;\;
f\in X^{\ast}\biggr\}. \]
Then for any $\da >0$ there exists continuous function
$l_{\da}\cn T\ra\mc{M}(X)$, the $\da -\text{\text{selection}}$ of
 $L$.
\et
\bw
(continuous version of the noncompact Krein -- Milman theorem)\\
Let $T$, $X$, $P$ and $p$ be as above. 
 Then for any $\g >0$ there exists continuous family of
measures $(\mu_{t})_{t\in T}$ on $X$, supported by the closure of the
set of extreme points of $P(t)$ and $\g -\text{\text{representing}}$
point $p(t)$.
\ew
The proofs of these results are identical to the ones given in the
previous section.

\tf{Acknowledgments.} Author would like to thank Andrzej
Fryszkowski, Piotr Mankiewicz and Kazimierz Br\k{a}giel for discussions.

\noindent Institute of Mathematics and Computer Science\\
Technical University of Cz\k{e}stochowa\\
J.~H.~D\k{a}browskiego 73\\
42-200 Cz\k{e}stochowa\\
Poland\\
E-mail: ppuchala@imi.pcz.pl

\end{document}